\documentclass[11pt,letterpaper,reqno]{amsart}
\usepackage{amsfonts}
\usepackage[normalem]{ulem}
\usepackage{xcolor}
\usepackage{bm} 
\usepackage{graphicx} 
\usepackage{tikz}

\usepackage{amsmath,amssymb,amsthm}
\usepackage{latexsym}
\usepackage[utf8x]{inputenc}
\usepackage[T1]{fontenc}
\usepackage{fancyhdr}

\newcommand{\R}{\mathbb{R}}

\usepackage[colorlinks=true,urlcolor=blue]{hyperref}

\newtheorem{theorem}{Theorem}[section]
\newtheorem{lemma}[theorem]{Lemma}
\newtheorem{proposition}[theorem]{Proposition}

\newcommand{\diam}{\operatorname{diam}}
\newcommand{\ip}[2]{\langle #1,#2\rangle}

\title[Distance between the centers of mass of a convex body and its boundary]{The sharp bound on the distance between the centers of mass of a convex body and its boundary}
\author{S. Myroshnychenko, D. Ryabogin, K. Tatarko, and V. Yaskin}

\address{S.~Myroshnychenko, Department of Mathematics and Statistics,  University of the Fraser Valley, Abbotsford, BC V2S 7M8,  Canada}

\email{serhii.myroshnychenko@ufv.ca}

\address{D.~Ryabogin, Department of Mathematics, Kent State University, Kent, OH 44242, USA}

\email{ryabogin@math.kent.edu}

\address{K.~Tatarko, Department of Pure Mathematics, University of Waterloo, Waterloo, ON N2L 3G1, Canada}

\email{ktatarko@uwaterloo.ca}

\address{V.~Yaskin, Department of Mathematical and Statistical Sciences, University of Alberta, Edmonton, AB T6G 2G1, Canada}

\email{yaskin@ualberta.ca}

\thanks{The first, third and fourth authors are supported in part by NSERC. The second  author is  supported in part by U.S.~National Science Foundation Grant DMS-2247771}



\begin{document}

	\begin{abstract}
		We prove that the distance between the centers of mass $c(K)$ and $c(\partial K)$ of a convex body $K$ and its boundary $\partial K$, respectively, does not exceed $ D_n$ times the width of $K$ in the direction of the vector joining $c(\partial K)$ and $c(K)$. The constant $D_n$ is optimal and coincides with the sharp constant in the problem of bounding the distance between the orthogonal projection of the centroid of a convex body onto a hyperplane and the centroid of its orthogonal projection onto that hyperplane.
	\end{abstract}
	
		\maketitle

	\section{Introduction}
Let $K\subset \mathbb R^n$ be a convex body, i.e., a convex  compact  set with non-empty interior. Write
\[
c(K)=\frac1{|K|}\int_K x\,dx,
\qquad
c(\partial K)=\frac1{|\partial K|}\int_{\partial K}x\,dS(x),
\]
for the centroid of $K$ and its boundary, respectively, where $|K|$ is the volume of $K$,  $|\partial K|$ is its surface area, and $dS$ denotes the surface area measure on $\partial K.$

 Let $S^{n-1}$ denote the unit sphere in $\mathbb{R}^n$. For $u\in S^{n-1}$, let
\[
w_K(u)=\max_{x\in K}\ip xu-\min_{x\in K}\ip xu
\]
be the width of $K$ in the direction $u$.

In their book \emph{Unsolved Problems in Geometry}, Croft, Falconer, and Guy \cite[p.~36, Problem~A25]{CFG} ask how far apart the centroid of a planar convex body and the centroid of its boundary can be, relative to such quantities as the diameter and the width.  They also point out the analogous questions in higher dimensions.

The planar problem was solved  by Nazarov, Ryabogin, and Yaskin \cite{NRY}.  They proved that for every planar convex body $K$ and every
$u\in S^1$,
\[
\bigl|\ip{c(\partial K)-c(K)}u\bigr|
\le \frac16\,w_K(u),
\]
and therefore
\[
|c(\partial K)-c(K)|
\le \frac16\,\diam K.
\]
The constant $1/6$ in both inequalities is best possible.

The purpose of the present paper is to solve the corresponding problem in all dimensions. In \cite{NRY}, it was conjectured that the supremum of
\[
\frac{|\langle c(\partial K) -  c(K),u\rangle|}{w_K(u)}
\]
over all convex bodies \(K\subset\mathbb R^n\) and all directions \(u\in S^{n-1}\) is the constant \(D_n\) defined below.

For an affine function $\ell(t) = (1-t)\alpha + t \beta$ on $[0,1]$ such that $\alpha, \beta \geq 0$ are not both zero, set
\begin{equation}\label{eq:Ij}
	I_j(\ell)=\int_0^1 t^j\ell(t)^{n-2}\,dt,\qquad j=0,1,2,
\end{equation}
and let
\begin{equation}\label{eq:Dn}
	D_n=
	\max_{\substack{\ell\ge0\ \mathrm{affine}\ \mathrm{on}\ [0,1]\\
			\ell\not\equiv0}}
	\left\{\frac{I_2(\ell)}{I_1(\ell)}-
	\frac{I_1(\ell)}{I_0(\ell)}\right\}.
\end{equation}
Multiplying \(\ell\) by a positive constant does not change the expression in \eqref{eq:Dn}. Hence we may normalize \(\ell\) so that \(\ell(0)+\ell(1)=1\). The resulting parameter set is compact, so the maximum is attained. The constants $D_n$ form an increasing sequence such that
\[
D_2=\frac16,\qquad
D_3=1-\sqrt{\frac23},\qquad \lim_{n\to \infty} D_n= 0.2016...,
\]
see \cite{MTY}.

The main result of this paper is the following.
\begin{theorem}\label{thm:main}
	For every convex body $K\subset\R^n$, $n\ge2$, and every $u\in S^{n-1}$,
	\begin{equation}\label{eq:width-main}
		\bigl|\ip{c(\partial K)-c(K)}u\bigr|\le D_n w_K(u).
	\end{equation}
	Consequently,
	\begin{equation}\label{eq:diam-main}
		|c(\partial K)-c(K)|\le D_n\diam K.
	\end{equation}
	The constant $D_n$ is optimal in both inequalities.
\end{theorem}

The constants $D_n$ first appeared in the following problem studied by Myroshnychenko, Tatarko, and Yaskin \cite{MTY}. For  a convex body $K\subset \mathbb R^n$ and  a hyperplane $H$ in $\mathbb R^n$, consider  $P_Hc(K)$, the orthogonal projection of $c(K)$ onto $H$, and $c(P_HK)$, the centroid of the orthogonal projection of $K$ onto $H$. If these two points are different, let 
$$
u = \frac{P_Hc(K)-c(P_HK)}{|P_Hc(K)-c(P_HK)|}.
$$
It was shown in \cite{MTY} that 
$$
|P_Hc(K)-c(P_HK)|\le D_n w_K(u).
$$
Here the constant $D_n$ is best possible.

There is a geometric reason to expect the constants $D_n$ from \cite{MTY} to appear in Theorem~\ref{thm:main}. Fix a hyperplane $H$ and compress a convex body $K$ in the direction orthogonal to $H$, preserving its projection onto $H$. The centroid of the compressed body projects to $P_Hc(K)$, while the projection of its boundary centroid  approaches $c(P_HK)$. Indeed, as the body becomes thinner, its surface area becomes concentrated on the two boundary graphs over $P_HK$. Thus the inequality of \cite{MTY} is recovered as a limiting case of the inequality in Theorem~\ref{thm:main}. In particular, the optimal constant in \eqref{eq:width-main} cannot be smaller than $D_n$.

The paper is organized as follows. Section~\ref{sec:outline} outlines the proof and states its main ingredient, Proposition~\ref{prop:base-comparison}. In Section~\ref{Sec:proof}, we derive bounds (\ref{eq:width-main}) and (\ref{eq:diam-main}) from this proposition, and in Section~\ref{sec:sharpness} we establish the sharpness of the constant $D_n$. Section~\ref{sec:two ii} contains one-dimensional auxiliary inequalities. Finally, Section~\ref{sec:prop} uses localization to prove Proposition~\ref{prop:base-comparison}.

\section{Outline of the proof}\label{sec:outline}

The proof of Theorem~\ref{thm:main} is based on decomposing $K$
into parallel chords. The coordinate of $c(K)$ in a fixed direction is expressed as an average of midpoint coordinates of chords weighted by chord length; see \eqref{eq:volume-chords}. To obtain the corresponding coordinate of $c(\partial K)$, we first take the unweighted average of these midpoint coordinates for each chord direction and then average the resulting values over all directions, with weights given by the volumes of the corresponding projections; see \eqref{eq:boundary-chords}. The following proposition allows us to compare these averages.
\begin{proposition}\label{prop:base-comparison}
	Let $\Omega\subset\R^{n-1}$ be a bounded, open, and convex set. Let $a,b$ be bounded
	 functions on $\Omega$, with $a$ convex, $b$ concave, and
	$0\le a\le b\le1$. Put
	\[
	h=b-a,\qquad m=\frac{a+b}{2}.
	\]
	If $h\not\equiv0$, then
	\begin{equation}\label{eq:base-comparison}
		\left|
		\frac1{|\Omega|}\int_\Omega m(y)\,dy-
		\frac{\int_\Omega h(y)m(y)\,dy}{\int_\Omega h(y)\,dy}
		\right|\le D_n.
	\end{equation}
\end{proposition}

The appearance of $D_n$ in Proposition~\ref{prop:base-comparison} can be understood through localization. Applied in $\mathbb R^{n-1}$, the localization lemma of Kannan,
Lov\'asz, and Simonovits reduces the proof of \eqref{eq:base-comparison} to one-dimensional integrals with weights proportional to $\ell^{n-2}$, where $\ell$ is nonnegative and affine. These are precisely the weights occurring in the definition of $D_n$. For affine functions, the required one-dimensional bound follows from this definition. The one-dimensional auxiliary inequalities extend this bound to general pairs of convex and concave functions.

\section{Proof of the bounds in Theorem~\ref{thm:main}}\label{Sec:proof}

\begin{proof}
	Assuming Proposition~\ref{prop:base-comparison}, we now prove inequalities \eqref{eq:width-main} and \eqref{eq:diam-main}.  Fix $u\in S^{n-1}$. Since the inequalities are invariant under translations and dilations, we may assume that
	\[
	\min_{x\in K}\ip xu=0,
	\qquad
	\max_{x\in K}\ip xu=1.
	\]
	Thus it suffices to prove
	\[
	|\ip{c(\partial K)-c(K)}u|\le D_n.
	\]
	
	Fix $v\in S^{n-1}$ with $\ip uv\ne0$. Replacing $v$ by $-v$ if necessary, we may assume $\ip uv>0$. Put
	\[
	\Omega_v=\operatorname{int}(P_{v^\perp}K),
	\]
	where $P_{v^\perp}$ denotes the orthogonal projection onto the hyperplane $v^\perp=\{x\in \mathbb R^n: \langle x,v\rangle =0\},$ and the interior is taken relative to $v^\perp$.
	
	For $y\in\Omega_v$, define
	\[
	a_v(y)=\min_{x\in K\cap(y+\R v)}\ip xu,
	\qquad
	b_v(y)=\max_{x\in K\cap(y+\R v)}\ip xu.
	\]
	By convexity of $K$, the function $a_v$ is convex and $b_v$ is concave on $\Omega_v$. Moreover,	$0\le a_v<b_v\le1$. Put
	\[
	h_v=b_v-a_v,
	\qquad
	m_v=\frac{a_v+b_v}{2}.
	\]
	Thus $m_v(y)$ is the $u$-coordinate of the midpoint of the chord $K\cap(y+\R v)$, while the length of this chord is $h_v(y)/\ip uv$.

    We first express $\ip{c(K)}u$ in terms of these functions. Writing $x=y+tv$, where $y\in v^\perp$, and using Fubini's theorem gives
	\[
	\begin{aligned}
		\int_K\ip xu\,dx
		&=\int_{\Omega_v}
		\left(
		\int_{\{t:y+tv\in K\}}
		\ip{y+tv}u\,dt
		\right)dy\\
		&=\frac1{\ip uv}
		\int_{\Omega_v}
		\left(\int_{a_v(y)}^{b_v(y)}s\,ds\right)dy\\
		&=\frac1{\ip uv}
		\int_{\Omega_v}
		\frac{b_v(y)^2-a_v(y)^2}{2}\,dy\\
		&=\frac1{\ip uv}
		\int_{\Omega_v}h_v(y)m_v(y)\,dy.
	\end{aligned}
	\]
Here, in the second equality, we used the change of variables $s = \ip{x}u$ so that $ds = \ip uv dt.$ Similarly,	
\[
	|K|=\frac1{\ip uv}\int_{\Omega_v}h_v(y)\,dy.
	\]
	Consequently,
	\begin{equation}\label{eq:volume-chords}
		\ip{c(K)}u
		=\frac1{|K|}\int_K\ip xu\,dx
		=
		\frac{\displaystyle\int_{\Omega_v}h_v(y)m_v(y)\,dy}
		{\displaystyle\int_{\Omega_v}h_v(y)\,dy}.
	\end{equation}

We next express $\ip{c(\partial K)}u$ using the same chord midpoint coordinates $m_v(y)$. Let $\nu_K(x)$ denote the outer unit normal to $K$, defined $dS$-almost everywhere on $\partial K$. For every nonnegative measurable function $f$, the area formula (see \cite[Theorem~3.2.3]{Federer}) gives
	\begin{equation}\label{eq:area-formula-projection}
		\int_{\partial K}
		f(x)|\ip{\nu_K(x)}v|\,dS(x)
		=
		\int_{v^\perp}
		\sum_{x\in\partial K\cap(y+\R v)}f(x)\,dy.
	\end{equation}
For $y\in\Omega_v$, the fiber on the right consists of the two endpoints of the corresponding chord. The boundary of $P_{v^\perp}K$ has $(n-1)$-dimensional Lebesgue measure zero and therefore does not contribute to the integral. We may apply \eqref{eq:area-formula-projection} with $f=1$ and $f(x)=\ip xu$, obtaining
	\begin{equation}\label{int1}
		\int_{\partial K}|\ip{\nu_K(x)}v|\,dS(x)
		=2|\Omega_v|
	\end{equation}
	and
	\begin{equation}\label{int2}
		\int_{\partial K}\ip xu\,|\ip{\nu_K(x)}v|\,dS(x)
		=2\int_{\Omega_v}m_v(y)\,dy.
	\end{equation}
	Since $v$ and $-v$ determine the same chords, the functions $a_v,b_v,h_v,m_v$, as well as $\Omega_v$, are unchanged under $v\mapsto -v$. Consequently, \eqref{eq:volume-chords}, \eqref{int1}, and \eqref{int2} hold for every $v$ with $\ip uv\ne0$.
	
	 To apply Proposition~\ref{prop:base-comparison}, we now average \eqref{int1} and \eqref{int2} over all chord directions $v$ to obtain a suitable representation of the boundary centroid. Let $d\sigma$ be the normalized rotation-invariant measure on $S^{n-1}$. Since
	\[
	\int_{S^{n-1}}|\ip uv|\,d\sigma(v)=\gamma_n>0
	\]
	is independent of the unit vector $u$, Fubini's theorem together with \eqref{int2} gives		
	\begin{align*}
		\int_{\partial K}\langle x, u\rangle \,dS(x)
		&	=  \frac1{\gamma_n} \int_{\partial K}\langle x, u\rangle 	\left(
		\int_{S^{n-1}}
		|\ip{\nu_K(x)}v|\,d\sigma(v)
		\right) \,dS(x)\\
		&	=  \frac1{\gamma_n} \int_{S^{n-1}}  	\left(
		\int_{\partial K} \langle x, u\rangle 
		|\ip{\nu_K(x)}v|\,dS(x)
		\right)  \,d\sigma(v)\\
		&	=  \frac2{\gamma_n} \int_{S^{n-1}}  	\left(
		\int_{\Omega_v}m_v(y)\,dy
		\right)  \,d\sigma(v).
	\end{align*}
	By \eqref{int1}, we have 	
	\begin{align*}\int_{S^{n-1}}|\Omega_v| \,d\sigma(v)&=\frac12\int_{S^{n-1}} \left(\int_{\partial K}|\ip{\nu_K(x)}v|\,dS(x)\right)  \,d\sigma(v) \\
		&	=  \frac12 \int_{\partial K} 	\left(
		\int_{S^{n-1}}
		|\ip{\nu_K(x)}v|\,d\sigma(v)
		\right) \,dS(x)=\frac{\gamma_n}{2} |\partial K|.
	\end{align*}
	Therefore,
	\begin{equation}\label{eq:boundary-chords}
		\begin{aligned}
			\ip{c(\partial K)}u
			=\frac1{|\partial K|}
			\int_{\partial K}\ip xu\,dS(x)=
			\frac{\displaystyle
				\int_{S^{n-1}}\int_{\Omega_v}
				m_v(y)\,dy\,d\sigma(v)}
			{\displaystyle
				\int_{S^{n-1}}|\Omega_v|\,d\sigma(v)}.
		\end{aligned}
	\end{equation}

	We now apply Proposition~\ref{prop:base-comparison}. Since $0\le a_v<b_v\le1$, it follows from \eqref{eq:volume-chords} that
	\[
	\begin{aligned}
		\left|
		\frac1{|\Omega_v|}\int_{\Omega_v}m_v(y)\,dy
		-\ip{c(K)}u
		\right|
		=
		\left|
		\frac1{|\Omega_v|}\int_{\Omega_v}m_v(y)\,dy
		-\frac{\displaystyle\int_{\Omega_v}
			h_v(y)m_v(y)\,dy}
		{\displaystyle\int_{\Omega_v}h_v(y)\,dy}
		\right|
		\le D_n
	\end{aligned}
	\]
	whenever $\ip uv\ne0$.
	
	Using \eqref{eq:boundary-chords}, we therefore obtain
	\[
	\begin{aligned}
		|\ip{c(\partial K)-c(K)}u|
		=
		\frac{\displaystyle
			\left|
			\int_{S^{n-1}}|\Omega_v|
			\left(
			\frac1{|\Omega_v|}\int_{\Omega_v}m_v(y)\,dy
			-\ip{c(K)}u
			\right)d\sigma(v)
			\right|}
		{\displaystyle\int_{S^{n-1}}|\Omega_v|\,d\sigma(v)}
		\le D_n.
	\end{aligned}
	\]
	Here the directions satisfying $\ip uv=0$ may be disregarded, since they form a set of spherical measure zero.

	The last inequality is precisely \eqref{eq:width-main} since $w_K(u)=1$. If $c(\partial K)\ne c(K)$, choosing $u$ in the
	direction of their difference and using $w_K(u)\le\operatorname{diam}K$ gives \eqref{eq:diam-main}.
\end{proof}

	\section{Sharpness}\label{sec:sharpness}

    To prove the sharpness of the constant $D_n$ in \eqref{eq:width-main}, it is enough to exhibit a sequence of bodies for which the quotient in \eqref{eq:width-main} tends to $D_n$ in the direction $e_1$.

	Assume first $n\ge 3$. Choose an affine function $\ell\ge0$ attaining the maximum in \eqref{eq:Dn}, and for $\varepsilon>0$ define
	\begin{equation}\label{eq:sharp-family}
		K_\varepsilon=
		\{(t,y,z)\in\R\times\R^{n-2}\times\R:
		0\le t\le1,\quad
		\|y\|_\infty\le\varepsilon\ell(t),\quad
		0\le z\le\varepsilon^2t\}.
	\end{equation}
	
	The body ${K_\varepsilon}$ is convex and $w_{K_\varepsilon}(e_1)=1$. Fubini's theorem gives
	\begin{equation}\label{eq:sharp-volume}
		\ip{c(K_\varepsilon)}{e_1}=\frac{I_2(\ell)}{I_1(\ell)}.
	\end{equation}
	 
    The bottom and top faces, $z=0$ and	$z=\varepsilon^2t$, have total area
	\[
	\bigl(1+\sqrt{1+\varepsilon^4}\bigr)(2\varepsilon)^{n-2} I_0(\ell),
	\]
	and  the first coordinates of their centroids are both $I_1(\ell)/I_0(\ell)$. The $2(n-2)$ side faces have total area $O(\varepsilon^{n-1})$. The part of the boundary at $t=1$ has area $O(\varepsilon^{n})$, and the part at $t=0$ has zero $(n-1)$-dimensional area. Hence
	\begin{equation}\label{eq:sharp-boundary}
		\ip{c(\partial K_\varepsilon)}{e_1}\longrightarrow
		\frac{I_1(\ell)}{I_0(\ell)}.
	\end{equation}

	Combining \eqref{eq:sharp-volume} and \eqref{eq:sharp-boundary} gives
	\[
	\frac{|\ip{c(\partial K_\varepsilon)-c(K_\varepsilon)}{e_1}|}
	{w_{K_\varepsilon}(e_1)}
	\longrightarrow
	\frac{I_2(\ell)}{I_1(\ell)}-\frac{I_1(\ell)}{I_0(\ell)}=D_n.
	\] 

    Thus the constant in \eqref{eq:width-main} is sharp for $n\ge 3$.
    
    Since $K_\varepsilon$ contains $(0,0,\dots,0)$ and $(1,0,\dots,0)$ and lies within distance $O(\varepsilon)$ of the segment joining these points, we have
	\[
	\diam K_\varepsilon\longrightarrow1.
	\]
	The $y$-coordinates of both centroids vanish by symmetry and their $z$-coordinates are $O(\varepsilon^2)$. Therefore
	\[
	\frac{|c(\partial K_\varepsilon)-c(K_\varepsilon)|}
	{\diam K_\varepsilon}
	\longrightarrow D_n,
	\]
	which proves sharpness in \eqref{eq:diam-main} for $n\ge 3$.
    
   When $n = 2$, let  $K_\varepsilon$ be the triangle with vertices $(0,0),(1,0),(1,\varepsilon^2)$. Then $\ip{c(K_\varepsilon)}{e_1}=2/3$, $\ip{c(\partial K_\varepsilon)}{e_1}\to 1/2$, and ${\diam K_\varepsilon}\to 1$ as $\varepsilon \to 0^+$. Thus we obtain the same conclusions as above.

    \qed
	
	\bigskip
	
	\section{One-dimensional comparison inequalities}\label{sec:two ii}
	
The remainder of the paper is devoted to proving Proposition~\ref{prop:base-comparison}. In this section we obtain a series of auxiliary one-dimensional integral inequalities that will be used in the proof of Proposition~\ref{prop:needle-comparison}. The latter will be the main ingredient in the  localization argument of Proposition~\ref{prop:base-comparison}.

\subsection{Preliminary estimates}
	
	For any nonzero finite positive measure $\nu$  and bounded measurable real-valued function $f$, define
	\[
	A_\nu(f)=\frac{\int f\,d\nu}{\int d\nu}
	\]
	and
	\[
V_\nu(f)=A_\nu(f^2)-A_\nu(f)^2.
	\]
	Thus $A_\nu$ denotes averaging, and $V_\nu$ is the	corresponding variance. Notice that
	\[
	V_\nu(f)
	=A_\nu\bigl((f-A_\nu(f))^2\bigr)\ge0.
	\]

	For the remainder of the paper, put $C=D_n$. Fix  a nonzero nonnegative affine function $\ell(t)$ on $[0,1]$ and define the normalized weight
	\begin{equation}\label{eq:weight}
		\omega(t)=\frac{\ell(t)^{n-2}}{I_0(\ell)},\qquad 0\le t\le1.
	\end{equation}
	Let $d\mu(t)=\omega(t)\,dt$. Throughout this section, we omit the subscript \(\mu\) and write \(A\) and \(V\)  for \(A_\mu\) and \(V_\mu\), respectively. Since  $\mu([0,1])=1$, we have
	\[
	A(f)=\int_0^1 f(t)\omega(t)\,dt.
	\]
	Let $T(t)=t$. Then
	\[
	A(T)=\int_0^1 t\,\omega(t)\,dt
	=\frac{I_1(\ell)}{I_0(\ell)},
	\]
	and
	\[
	A(T^2)=\int_0^1 t^2\,\omega(t)\,dt=\frac{I_2(\ell)}{I_0(\ell)}.
	\]
	Therefore
	\[
	V(T)
	=A(T^2)-A(T)^2
	=
	\frac{I_2}{I_0}-\frac{I_1^2}{I_0^2}.
	\]
	From the definition of $D_n$,
	\[
	\frac{V(T)}{A(T)}
	=
	\frac{I_2(\ell)}{I_1(\ell)}
	-
	\frac{I_1(\ell)}{I_0(\ell)}
	\le C,
	\]
	and hence
	\[
	V(T)\le CA(T).
	\]
	Since $A$ is linear and $V(1-T) = V(T)$, applying the same inequality after the reflection $t\mapsto1-t$ gives
	\[
	V(T)\le C\bigl(1-A(T)\bigr).
	\]
	Thus
	\begin{equation}\label{eq:moment-bounds}
		V(T)\le CA(T),\qquad\mbox{and}\qquad
		V(T)\le C\bigl(1-A(T)\bigr).
	\end{equation}
	We shall also use
	\[
	C\ge\frac16,
	\]
	which follows by taking $\ell$ constant in \eqref{eq:Dn}.

In the lemmas below, we establish integral inequalities for bounded convex or concave functions on $[0,1]$. Since changing their endpoint values does not affect these integrals, we may assume that such functions are continuous on $[0,1]$, replacing their endpoint values by their one-sided limits if necessary.

The following estimate will be used in Lemma~\ref{lem:paired-hinges}.

\begin{lemma}\label{lem:concave-weight}
	If $0\le h\le1$ is a concave function on $[0,1]$, then
	\begin{equation}\label{eq:concave-variance}
		V(h)\le\frac13 A(h).
	\end{equation}
\end{lemma}

\begin{proof}	 Suppose first that $n\ge 3$, and consider
	the convex body
	\[
	P_\ell=\{(t,y)\in[0,1]\times\R^{n-2}:
	\|y\|_\infty\le\ell(t)\}.
	\]
	For fixed $t$, the $y$-section of $P_\ell$ has $(n-2)$-dimensional volume
\(
(2\ell(t))^{n-2}.
\)
Note that the explicit construction of $P_\ell$ is not important. All we need is that the $(n-2)$-dimensional volume of its cross-sections is proportional to $\ell^{n-2}(t)$.

	Hence, for every integrable function $f$ on $[0,1]$,
	\[
	\frac1{|P_\ell|}\int_{P_\ell}f(t)\,dt\,dy
	=
	\frac{\int_0^1 f(t)\ell(t)^{n-2}\,dt}
	{\int_0^1\ell(t)^{n-2}\,dt}
	=A(f).
	\]
	When $n=2$, we take $P_\ell=[0,1]$, and the same identity holds.
	
	If $A(h)=0$ or $A(h)=1$, then $V(h)=0$ and there is nothing to prove. Hence for the remainder of the proof we assume
	$0<A(h)<1$.
	
	Let
	\[
	F(s)=A\bigl(\mathbf 1_{\{h>s\}}\bigr)=\frac{|\{(t,y)\in P_\ell:h(t)>s\}|}{|P_\ell|},
	\qquad 0\le s\le1.
	\]	
	Since $h$ is nonnegative, concave, and not identically
	zero, it is positive on $(0,1)$, so $F(0)=1$.

	Since $(t,y)\mapsto h(t)$ is concave on $P_\ell$, the Brunn--Minkowski inequality implies that $F$ is log-concave.
	By the layer cake decomposition,
	\begin{equation}\label{eq:layer-cake-concave}
		A(h)=\int_0^1F(s)\,ds,
		\qquad\mbox{and}\qquad
		A(h^2)=2\int_0^1sF(s)\,ds.
	\end{equation}
	 To bound $V(h),$ we compare the function $F$ with an exponential function that has the same integral. In particular, this will allow us to bound $A(h^2)$ in terms of $A(h)$. Choose $\lambda>0$ so that
	\[
	A(h)=\int_0^1e^{-\lambda s}\,ds
	=\frac{1-e^{-\lambda}}{\lambda},
	\]
	which is always possible since $\frac{1-e^{-\lambda}}{\lambda}$  decreases continuously from 1 to 0 as $\lambda$ increases from $0$ to $\infty$.

	Since $\log F(s)+\lambda s$ is concave and vanishes at $s=0$, the difference $F(s)-e^{-\lambda s}$ changes sign at most once, from nonnegative to nonpositive. Thus there exists
	$s_0\in[0,1]$ such that
	\[
	(s-s_0)\bigl(F(s)-e^{-\lambda s}\bigr)\le0
	\]
	for all $s\in[0,1]$. Since
	\[
	\int_0^1F(s)\,ds=\int_0^1e^{-\lambda s}\,ds,
	\]
	we obtain
	\[
		\int_0^1sF(s)\,ds-\int_0^1se^{-\lambda s}\,ds
		=\int_0^1(s-s_0)
		\bigl(F(s)-e^{-\lambda s}\bigr)\,ds
		\le0.
	\]
Together with \eqref{eq:layer-cake-concave}, this implies
	\[
	A(h^2)
	=2\int_0^1sF(s)\,ds
	\le 2\int_0^1se^{-\lambda s}\,ds
	=\frac{2}{\lambda^2}
	\bigl(1-(1+\lambda)e^{-\lambda}\bigr).
	\]

Therefore,
\begin{equation}\label{eq:lambda-bound}
\frac{V(h)}{A(h)}
=\frac{A(h^2)}{A(h)}-A(h)\le \frac{2 \bigl(1-(1+\lambda)e^{-\lambda}\bigr)}{\lambda (1-e^{-\lambda})}
-\frac{1-e^{-\lambda}}{\lambda}
\le\frac13.
\end{equation}
	For the last inequality, multiplication by the positive denominator reduces it to
	\[
	\Phi(\lambda)=(\lambda-3)e^\lambda+3e^{-\lambda}+5\lambda\ge0.
	\]
	Now $\Phi(0)=\Phi'(0)=0$ and
	\[
	\Phi''(\lambda)=(\lambda-1)e^\lambda+3e^{-\lambda}>0.
	\]
	For $\lambda\ge1$ this is immediate. For $0\le\lambda\le1$, it follows from$(1-\lambda)e^{2\lambda}\le e/2<3$. Thus $\Phi$ is increasing and so \eqref{eq:lambda-bound} holds.

\end{proof}

	The goal of this section is to establish  Proposition \ref{prop:needle-comparison}, i.e., we need to prove the inequality
	\[
-2CA(b-a)\le	V(a)-V(b)\le 2CA(b-a)
	\] 
	for continuous functions $0\le a\le b\le1$ on $[0,1]$, with $a$ convex and $b$ concave. Consider the functional 
	\[
	J(f)=V(f)+2CA(f).
	\] 	
	Then the right-hand inequality above is equivalent to $$J(a)\le J(b).$$ The advantage of this formulation is that $J$ satisfies certain comparison inequalities for elementary piecewise affine functions, which can then be extended to general convex functions using the convexity of $J$.

	To see that the functional $J$ is convex, take $0\le\lambda\le1$. Then
	\[
	V(\lambda f+(1-\lambda)g)
	=
	\lambda V(f)+(1-\lambda)V(g)
	-\lambda(1-\lambda)V(f-g)
	\le
	\lambda V(f)+(1-\lambda)V(g),
	\]
	while $A$ is linear. Hence
	\[
	J(\lambda f+(1-\lambda)g)
	\le
	\lambda J(f)+(1-\lambda)J(g).
	\]
	
	\subsection{One-sided elementary piecewise affine functions}
	
	For $0\le s<1$, set
	\[
	R_s(t)=\frac{(t-s)_+}{1-s},
	\]
	and let $R_1=0$.

	\begin{lemma}\label{lem:hinge}
		The function $s\mapsto J(R_s)$ is nonincreasing on $[0,1]$. In particular,
		\begin{equation}\label{eq:hinge}
			J(R_s)\le J(R_0)=J(T).
		\end{equation}
	\end{lemma}
	
	\begin{proof}
   
     We first show that $J(R_s)$ is nonincreasing on $s \in [0,1)$ by proving that $\frac d{ds}J(R_s)\le0$.  For $0 \leq s<1$, let $\mu_s$ be the probability measure on $[s,1]$ defined by
    \[
		d\mu_s(t)=
		\frac{\omega(t)\,dt}{\displaystyle\int_s^1\omega(r)\,dr}.
		\] 
        Also, for $s \leq t\leq1$,
		\[
		R_s(t)=\frac{t-s}{1-s}.
		\]
		Under the change of variables
		\[
		r=\frac{t-s}{1-s},
		\]
        \(R_s(t)\) becomes the coordinate function \(T(r)=r\) on \([0,1]\), and the measure $\mu_s$ is now defined on $[0,1]$ and is given by
		\[
	\frac{\ell(s+(1-s)r)^{\,n-2}}
		{\displaystyle\int_0^1
		\ell(s+(1-s)\rho)^{\,n-2}\,d\rho}\,dr.
		\]
        Since $r\mapsto\ell(s+(1-s)r)$ is again a nonnegative affine function on $[0,1]$, this is a weight of the form \eqref{eq:weight}. Hence the second inequality in \eqref{eq:moment-bounds}, applied to $\mu_s$, gives
		\begin{equation}\label{eq:restricted-bound}
			V_{\mu_s}(R_s)
			\le
			C\bigl(1-A_{\mu_s}(R_s)\bigr).
		\end{equation}

		We now differentiate $J(R_s)$. Recall that
		\[
		A(R_s)=\int_s^1 R_s(t)\omega(t)\,dt,
		\qquad
		A(R_s^2)=\int_s^1 R_s(t)^2\omega(t)\,dt.
		\]
		Observe that for $s< t\le1$,
		\[
		\frac{\partial R_s(t)}{\partial s}
		=-\frac{1-R_s(t)}{1-s}.
		\]
		By Leibniz's rule,
		\[
		\begin{aligned}
			\frac d{ds}A(R_s)
			&=-R_s(s)\omega(s)
			+\int_s^1\frac{\partial R_s(t)}{\partial s}\omega(t)\,dt\\
			&=-\frac1{1-s}\int_s^1(1-R_s(t))\omega(t)\,dt\\
			&=-\frac{\mu([s,1])}{1-s}
			\bigl(1-A_{\mu_s}(R_s)\bigr),
		\end{aligned}
		\]
		where we used $R_s(s)=0$  and $d\mu(t)=\omega(t) \,dt$.  Similarly,
		\[
			\frac d{ds}A(R_s^2) =
			\int_s^1 2R_s(t)\frac{\partial R_s(t)}{\partial s}\omega(t)\,dt
			=
			-\frac{2\mu([s,1])}{1-s}
			\left(
			A_{\mu_s}(R_s)-A_{\mu_s}(R_s^2)
			\right).
		\]
		
		Since
		\[
		V(R_s)=A(R_s^2)-A(R_s)^2
		\]
		and
		\[
		A(R_s)=\mu([s,1])A_{\mu_s}(R_s),
		\]
		substitution of the preceding two derivatives and basic rearrangements give
		\[
		\begin{aligned}
			\frac d{ds}J(R_s)
			=-\frac{2\mu([s,1])}{1-s}
			\Big[
			\bigl(C+(1-\mu([s,1]))A_{\mu_s}(R_s)\bigr)
			\bigl(1-A_{\mu_s}(R_s)\bigr)
			-V_{\mu_s}(R_s)
			\Big].
		\end{aligned}
		\]
		The expression in brackets can be written as
		\[
		\Big[
		C\bigl(1-A_{\mu_s}(R_s)\bigr)
		-V_{\mu_s}(R_s)
		\Big]
		+(1-\mu([s,1]))A_{\mu_s}(R_s)
		\bigl(1-A_{\mu_s}(R_s)\bigr).
		\]
		Both terms are nonnegative: the first by
		\eqref{eq:restricted-bound}, and the second because
		\[
		0\le\mu([s,1])\le1,
		\qquad
		0\le A_{\mu_s}(R_s)\le1.
		\]
		Therefore
		\[
		\frac d{ds}J(R_s)\le0.
		\]
 Thus $s\mapsto J(R_s)$ is nonincreasing on $[0,1)$. Since $J(R_s) \geq 0 = J(R_1)$,  \eqref{eq:hinge} holds on $[0,1].$ This completes the proof.
 
	\end{proof}

	\subsection{Two-sided elementary piecewise affine functions}
	
	For $0<u\le1$, put
	\[
	L_u(t)=\frac{(u-t)_+}{u},
	\]
	and set $L_0=0$.
	
	\begin{lemma}\label{lem:paired-hinges}
		If $0\le u\le s\le1$ and $0\le z\le1$, then
		\begin{equation}\label{eq:paired-hinges}
			J(zL_u+R_s)\le J\bigl(z(1-T)+T\bigr).
		\end{equation}
	\end{lemma}

	\begin{proof}
		Since $u\le s$, we have $L_uR_s=0$. Expanding the difference of the two sides of \eqref{eq:paired-hinges} as a quadratic polynomial in $z$ and regrouping gives
		\[
		\begin{aligned}
			&J\bigl(z(1-T)+T\bigr)-J(zL_u+R_s)\\
			&\quad=(1-z)^2\bigl[J(T)-J(R_s)\bigr]\\
			&\qquad+2z(1-z)\Bigl[
			C\bigl(1+A(T)-A(L_u)-2A(R_s)\bigr)
			-V(R_s)+A(L_u)A(R_s)
			\Bigr]\\
			&\qquad+z^2\bigl[J(1)-J(L_u+R_s)\bigr].
		\end{aligned}
		\]
		To show \eqref{eq:paired-hinges}, it is enough to show that every coefficient of $(1-z)^2, \ 2z(1-z), \ z^2$ is nonnegative. The first coefficient is nonnegative by Lemma~\ref{lem:hinge}.
		
		For the last coefficient, put
		\[
		h=1-L_u-R_s.
		\]
		Then $0\le h\le1$ is concave,
		\[
		V(L_u+R_s)= V(1-h) = V(h),
		\qquad
		A(h)=1-A(L_u)-A(R_s),
		\]
		and hence
		\[
		J(1)-J(L_u+R_s)
		=2CA(h)-V(h)\ge0
		\]
		by Lemma~\ref{lem:concave-weight} and $2C\ge1/3$.

	It remains to show that the mixed coefficient is nonnegative.
	Put $q=A(L_u)$. Since $0\le L_u+R_s\le1$, we have
	\[
	0\le q\le1-A(R_s).
	\]
	For fixed $s$, the mixed coefficient
	\[
	C\bigl(1+A(T)-q-2A(R_s)\bigr)
	-V(R_s)+qA(R_s)
	\]
	is an affine function of  $q$ on $[0,1-A(R_s)]$. It therefore suffices to check its values at the endpoints of this interval.
	
	At $q=0$, it equals
	\[
	\begin{aligned}
		C\bigl(1+A(T)-2A(R_s)\bigr)-V(R_s)
		=
		\bigl[J(T)-J(R_s)\bigr]
		+\bigl[C(1-A(T))-V(T)\bigr]\ge0
	\end{aligned}
	\]
	by Lemma~\ref{lem:hinge} and \eqref{eq:moment-bounds}.
	
	At $q=1-A(R_s)$, it equals
	\[
	C\bigl(A(T)-A(R_s)\bigr)
	+A(R_s)\bigl(1-A(R_s)\bigr)-V(R_s)\ge0,
	\]
	since $R_s\le T$ and
\[
V(R_s)=A(R_s^2)-A(R_s)^2\le A(R_s)-A(R_s)^2= A(R_s)\bigl(1-A(R_s)\bigr).
\]
	Thus all three coefficients are nonnegative, proving the claim.

	\end{proof}

	\subsection{Affine majorants}
	
	\begin{lemma}\label{lem:affine-majorant}
		Let $b$ be affine and let $a$ be  convex on $[0,1]$. If
		$0\le a\le b\le1$, then
		\begin{equation}\label{eq:affine-majorant}
			J(a)\le J(b).
		\end{equation}
	\end{lemma}

	\begin{proof}
    
   The idea of the proof is to approximate $a$ by convex piecewise affine functions, and then use convexity of $J$ and Lemma~\ref{lem:paired-hinges}.

    We first reduce to the case where $a$ vanishes at some point.
Subtracting $\min_{[0,1]}a$ from both $a$ and $b$ preserves
all the assumptions of the lemma and leaves $J(b)-J(a)$ unchanged,
since $V$ is invariant under addition of constants and
$A$ is linear. Thus we may assume that
\[
\min_{[0,1]}a=0.
\]

If $b\equiv0$, then $a\equiv0$ and there is nothing to prove.
Reflecting both the functions and the weight if necessary, we may assume $b(1)\ge b(0)$. We may further assume that $b(1)=1$. Indeed, if \eqref{eq:affine-majorant} holds for the functions $a/b(1)$ and $b/b(1)$, then
\[
V(a)-V(b)
\le 2C b(1)A(b-a)
\le 2CA(b-a),
\]
which is equivalent to $J(a)\le J(b)$. Here we used the linearity of $A$, the fact that $V$ is homogeneous of degree two, and the inequalities $0<b(1)\le1$ and $A(b-a)\ge0$.

These reductions preserve the condition $\min_{[0,1]}a=0$.
We therefore have
\[
b(t)=b(0)(1-t)+t,\qquad 0\le b(0)\le1,
\]
and $a(t_0)=0$ for some $t_0\in[0,1]$.

		For each $N$, subdivide $[0,t_0]$ and $[t_0,1]$ into $N$ equal subintervals:
		\[
		0=u_0<\cdots<u_N=t_0=s_0<\cdots<s_N=1
		\]
		(with the degenerate side omitted if $t_0=0$ or $1$), and let $p_N$ be the
		piecewise affine function equal to $ a$ at all partition points. Then $p_N$ is convex,
		\[
		p_N(t_0)=0,\qquad 0\le p_N\le b,
		\]
		and $p_N\to a$ uniformly.
		
		Write
		\[
		p_N=p_N^-+p_N^+,
		\]
		where
		\[
		p_N^-(t)=
		\begin{cases}
			p_N(t), & 0\le t\le t_0,\\
			0, & t_0\le t\le1,
		\end{cases}
		\qquad
		p_N^+(t)=
		\begin{cases}
			0, & 0\le t\le t_0,\\
			p_N(t), & t_0\le t\le1.
		\end{cases}
		\]
		Since $p_N\ge0$ and $p_N(t_0)=0$, both functions are continuous and convex.
		
		Moreover, $p_N^-$ vanishes on $[t_0,1]$, so it is a linear combination with nonnegative coefficients of the functions $(u_j-t)_+$, where $u_j\le t_0$. The coefficients are the jumps in its slope. Similarly, $p_N^+$ is a linear combination with nonnegative coefficients of the functions $(t-s_i)_+$, where $s_i\ge t_0$. Thus we may write
		\[
		p_N^-=\sum_j\alpha_j b(0)L_{u_j},
		\qquad
		\alpha_j\ge0,\quad \sum_j\alpha_j\le1,\quad u_j\le t_0,
		\]
		and
		\[
		p_N^+=\sum_i\beta_iR_{s_i},
		\qquad
		\beta_i\ge0,\quad \sum_i\beta_i\le1,\quad s_i\ge t_0.
		\]

		Here $\sum_i\beta_i\le1$ follows from $p_N^+(1)\le1$. If $b(0)>0$, then $\sum_j\alpha_j\le1$ follows from
		$p_N^-(0)\le b(0)$. If $b(0)=0$, then $t_0=0$ and $p_N^-\equiv0$; in this case we take a single term with $\alpha_0=1$ and $u_0=0$.

      In order to apply Lemma~\ref{lem:paired-hinges} and use convexity of $J$, we now express $p_N$ as a convex combination. Since $L_0=0$ and $R_1=0$, we may add them to the preceding sums so that
		\[
		\sum_j\alpha_j=\sum_i\beta_i=1.
		\]
		Hence
		\[
		p_N=\sum_{i,j}\alpha_j\beta_i
		\bigl(b(0)L_{u_j}+R_{s_i}\bigr).
		\]
		By Lemma~\ref{lem:paired-hinges} and the convexity of $J$,
		\[
		J(p_N)
		\le\sum_{i,j}\alpha_j\beta_i
		J\bigl(b(0)L_{u_j}+R_{s_i}\bigr)
		\le\sum_{i,j}\alpha_j\beta_iJ(b)=J(b).
		\]
		Letting $N\to\infty$ gives  
		$J(a)\le J(b)$.
	\end{proof}

	\begin{proposition}\label{prop:needle-comparison}
		Let $a$ be convex and $b$ be concave on $[0,1]$ with
		$0 \leq a\le b \leq 1$. Then
		\begin{equation}\label{eq:needle-comparison}
			|V(a)-V(b)|\le2CA(b-a).
		\end{equation}
	\end{proposition}

	\begin{proof}
		We first prove that $J(a)\le J(b)$. It suffices to do so when $b>0$ on $[0,1]$: the general case follows by applying the result to $(a+\delta)/(1+\delta)$ and $(b+\delta)/(1+\delta)$ and letting $\delta\to 0^+$, using the continuity of $J$ under uniform convergence.
		
		Let $L(t)=(1-t)b(0)+tb(1)$. Then
		\[
		0\le a\le L\le b\le1,
		\]
		and Lemma~\ref{lem:affine-majorant} gives $J(a)\le J(L)$.
		Since $L>0$ on $[0,1]$, for sufficiently small $\varepsilon>0$ the function
		\[
		g=(1+\varepsilon)L-\varepsilon b
		\]
		is nonnegative. It is also convex and satisfies $g\le L$, so Lemma~\ref{lem:affine-majorant} gives $J(g)\le J(L)$. Since
		\[
		L=\frac{g+\varepsilon b}{1+\varepsilon},
		\]
		the convexity of $J$ yields
		\[
		J(L)
		\le\frac{J(g)+\varepsilon J(b)}{1+\varepsilon}
		\le\frac{J(L)+\varepsilon J(b)}{1+\varepsilon}.
		\]
		Thus $J(L)\le J(b)$, and hence $J(a)\le J(b).$ Equivalently,
		\begin{equation}\label{eq:L4.5}
		    V(a)-V(b)\le 2CA(b-a).
		\end{equation}
	Applying the same argument to the convex function $1-b$ and the concave function $1-a$, satisfying $0\le1-b\le1-a\le1$, and using
\[
V(1-f)=V(f),
\qquad
(1-a)-(1-b)=b-a,
\]
gives
\[
V(b)-V(a)\le2CA(b-a).
\]
Together with \eqref{eq:L4.5}, this proves \eqref{eq:needle-comparison}.
	\end{proof}

	\section{Proof of 	
	Proposition~\ref{prop:base-comparison}}\label{sec:prop}

	Since $h\ge0$ is concave and not identically zero, $h>0$ in $\Omega$. Throughout the proof, $A$ and $V$ denote the functionals associated with Lebesgue measure on $\Omega$. Thus,
	\[
	A(f)=\frac1{|\Omega|}\int_\Omega f(y)\,dy,
	\qquad
	V(f)=A(f^2)-A(f)^2.
	\]
	Observe that the quantity we need to bound is
		\[
		\frac1{|\Omega|}\int_\Omega m\,dy-
		\frac{\int_\Omega hm\,dy}{\int_\Omega h\,dy}
		=
		A(m)-\frac{A(hm)}{A(h)}.
		\]
	Since $A$ is linear and
		\[
		h=b-a,\qquad m=\frac{a+b}{2},
		\]
	we obtain
		\[
		A(m)=\frac{A(a)+A(b)}2,
		\qquad
		hm=\frac{b^2-a^2}{2},
		\qquad
		A(h)=A(b)-A(a).
		\]
		Consequently,
		\begin{equation}\label{eq:midpoint-Q-identity}
			A(m)-\frac{A(hm)}{A(h)}
			=\frac12\left(
			A(a)+A(b)
			-
			\frac{A(b^2)-A(a^2)}
			{A(b)-A(a)}
			\right) =
			\frac{V(a)-V(b)}{2A(h)}.
		\end{equation}
		
	 We consider the planar case and the higher-dimensional case separately. If $n=2$, then $\Omega\subset\mathbb R$ is an interval. After an affine reparametrization of $\Omega$ by $(0,1)$, normalized Lebesgue measure on $\Omega$ becomes the constant weight, which is precisely the weight $\ell^{\,n-2}\equiv1$ from Section~\ref{sec:two ii}. Extending  \(a\) and \(b\) continuously to $[0,1]$ and applying Proposition~\ref{prop:needle-comparison}, we obtain
		\[
		|V(a)-V(b)|\le 2CA(h).
		\]
		Combining this with \eqref{eq:midpoint-Q-identity}, we obtain
		\[
		\left|
		\frac1{|\Omega|}\int_\Omega m\,dy-
		\frac{\int_\Omega hm\,dy}{\int_\Omega h\,dy}
		\right|=\left|
		A(m)-\frac{A(hm)}{A(h)}
		\right|
		\le C.
		\]
		This proves \eqref{eq:base-comparison} when $n=2$.

		We may therefore assume $n\ge3$ in the localization argument below. To reach a contradiction, suppose  that 
        $$ 
        A(m)>\frac{A(hm)}{A(h)}+C.
        $$ 
        Then the continuous functions
		\[
		\varphi=m-\frac{A(hm)}{A(h)}-C,\qquad \psi=h\left(m-\frac{A(hm)}{A(h)}\right)
		\]
		satisfy
		\[
		\int_\Omega\varphi>0,
		\qquad
		\int_\Omega\psi=0.
		\]

By the localization lemma of Kannan, Lov\'asz, and Simonovits \cite[Corollary~2.4]{KLS}, there exist a line segment $I\subset\Omega$ and a nonnegative affine function $\ell$ on $I$ such that
\begin{equation}\label{eq:Sec5two_int}
\int_I\varphi(x)\ell^{n-2}(x)\,dx>0,
\qquad
\int_I\psi(x)\ell^{n-2}(x)\,dx=0,
\end{equation}
where $dx$ denotes length measure along $I$.

Let $\nu$ denote the normalized measure on $I$ with density proportional to $\ell^{n-2}$. Using \eqref{eq:Sec5two_int} together with the definitions of $\varphi$ and $\psi$, we get
\[
A_\nu(m)>\frac{A(hm)}{A(h)}+C,
\qquad
\frac{A_\nu(hm)}{A_\nu(h)}
=\frac{A(hm)}{A(h)}.
\]
Here $A_\nu(h)>0$ because $h>0$ on $\Omega$. After an affine parametrization of $I$ by $[0,1]$, the measure $\nu$ has a density of the form \eqref{eq:weight}. The restrictions of $a$ and $b$ remain respectively convex and concave and satisfy
$0\le a\le b\le1$. Thus Proposition~\ref{prop:needle-comparison} applies.
	
The same algebra as in \eqref{eq:midpoint-Q-identity}, now applied to $\nu$, together with Proposition~\ref{prop:needle-comparison}, yields
		\[
		\left|
		A_\nu(m)-\frac{A_\nu(hm)}{A_\nu(h)}
		\right|
		=
		\frac{|V_\nu(a)-V_\nu(b)|}{2A_\nu(h)}
		\le C,
		\]
		a contradiction.

		To obtain the opposite inequality, apply the same argument to the functions
		\[
		\tilde a=1-b,\qquad \tilde b=1-a,
		\]
		where $\tilde a$ is convex, $\tilde b$ is concave, and $0\le \tilde a\le \tilde b\le 1$.
		For this pair,
		\[
		\tilde h= \tilde b - \tilde a=h,\qquad
		\tilde m=\frac{\tilde a+\tilde b}2=1-m.
		\]
		The inequality already proved therefore gives
		$$ 
		A(\tilde m)\le\frac{A(\tilde h \tilde m)}{A(\tilde h)}+C.
		$$
		Since
		$$
		A(\tilde m)=1 -A(m),\qquad\mbox{and}\qquad \frac{A(\tilde h \tilde m)}{A(\tilde h)}=1-\frac{A(hm)}{A(h)},
		$$
		we obtain
		\[
		A(1-m)\le1-\frac{A(hm)}{A(h)}+C,
		\]
		or equivalently
		\[
		A(m)\ge \frac{A(hm)}{A(h)}-C.
		\]
		This proves the opposite inequality and hence
		\eqref{eq:base-comparison}.
		\qed
 \bigskip

\noindent{\bf Acknowledgments.} The project was made possible by SQuaRE at the American Institute of Mathematics. The authors thank AIM for providing a supportive and mathematically rich environment.

    \medskip
	
	\noindent{\bf AI Disclaimer.} ChatGPT was used in exploratory discussions and in the preparation of this manuscript.

\end{document}